\documentclass{amsart}

\usepackage{amsmath, amssymb, amsthm}
\usepackage{tikz-cd}
\usepackage[shortlabels]{enumitem}
\usepackage[hidelinks]{hyperref}

\newtheorem{theorem}{Theorem}[section]
\newtheorem{lemma}[theorem]{Lemma}
\newtheorem{corollary}[theorem]{Corollary}
\newtheorem{proposition}[theorem]{Proposition}
\newtheorem{remark}[theorem]{Remark}

\newtheorem{fact}[theorem]{Fact}

\newtheorem{Assumption}[theorem]{Assumption}

\newtheorem*{theorem*}{Theorem}
\newtheorem*{corollary*}{Corollary}

\theoremstyle{definition}
\newtheorem{example}[theorem]{Example}
\newtheorem{definition}[theorem]{Definition}
\newtheorem*{definition*}{Definition}

\newtheorem{innercustomthm}{Theorem}
\newenvironment{customthm}[1]
  {\renewcommand\theinnercustomthm{#1}\innercustomthm}
  {\endinnercustomthm}

\newcommand{\calL}{\mathcal{L}}

\newcommand{\calC}{\mathcal{C}}

\newcommand{\calP}{\mathcal{P}}

\newcommand{\calQ}{\mathcal{Q}}
\newcommand{\calR}{\mathcal{R}}
\newcommand{\calS}{\mathcal{S}}

\newcommand{\calU}{\mathcal{U}}

\newcommand{\Cb}{\operatorname{Cb}}

\newcommand{\Aut}{\mathrm{Aut}}
\def\acl{\operatorname{acl}}
\def\Fib{\operatorname{Fib}}
\def\dcl{\operatorname{dcl}}
\def\tp{\operatorname{tp}}
\def\stp{\operatorname{stp}}

\newcommand{\fibre}[2]{\Fib_{f({#2})}(#1)}

\newcommand{\domequiv}{\hspace{0.7mm}{\scriptstyle\underline{\square}}\hspace{0.7mm}}
\def\forkindep{\mathrel{\raise0.2ex\hbox{\ooalign{\hidewidth$\vert$\hidewidth\cr\raise-0.9ex\hbox{$\smile$}}}}}

\newcommand{\partials}[2]{\frac{\partial #1}{\partial #2}}

\title{Domination, fibrations and splitting}

\author{Christine Eagles}
\author{L\'eo Jimenez}


\date{\today}

\begin{document}

\begin{abstract}
    This article is concerned with finite rank stability theory, and more precisely two classical ways to decompose a type using minimal types. The first is its domination equivalence to a Morley product of minimal types, and the second its semi-minimal analysis, both of which are useful in applications. Our main interest is to explore how these two decompositions are connected. We prove that neither determine the other in general, and give more precise connections using various notions from the model theory literature such as uniform internality, proper fibrations and disintegratedness. 
\end{abstract}

\maketitle

\tableofcontents

\section{Introduction}

In finite rank stability theory, there are two ways to decompose a type using minimal types, i.e. stationary types with only algebraic forking extensions. The first one, domination-decomposition, makes any finite Lascar rank type domination equivalent to a Morley product of minimal types. Recall that two stationary types $p,q$ are domination equivalent if there is some $D$ containing their parameters, some $a \models p\vert_D$ and $b \models q\vert_D$ such that for all $e$, we have $e \forkindep_D a$ if and only if $e \forkindep_D b$. This shows in particular that orthogonality is completely controlled by minimal types, an observation that was crucial in \cite{freitag2025finite} to produce non-orthogonality bounds for solutions of algebraic ordinary differential equations. Recent work on domination has been focused mostly on the unstable context, see for example \cite{haskell2008stable} and \cite{mennuni2022domination}, but this article is only considering superstable theories and types of finite Lascar rank. 

The second decomposition is the semi-minimal analysis. Given two types $p,q \in S(A)$, by a definable map $f : p \rightarrow q$, we mean an $A$-definable function with domain a formula in $p$ and codomain a formula in $q$. We may sometimes write $q$ as $f(p)$. The \emph{fibers} of $f$ are the complete types $\tp(b/f(b)A)$, for $b$ a realization of $p$. Recall that a stationary type $p \in S(A)$ is almost internal to some $r \in S(B)$ if there are $D \supset A \cup B$, some $a \models p \vert_D$ and some $c_1, \cdots , c_n \models r$ such that $a \in \acl(D, c_1, \cdots , c_n)$. A type is semi-minimal if it is almost internal to a minimal type. The semi-minimal analysis shows that these are the building block of finite rank types, as it gives a sequence of types and definable maps:
\[p=p_n\xrightarrow{f_{n-1}}p_{n-1}\xrightarrow{f_{n-2}}\cdots\xrightarrow{f_1}p_1\]
\noindent where for each $b_i \models p_i$, the strong type $\stp(b_i/f_{i-1}(b_i)A)$, as well as $p_1$, are semi-minimal.  Studying the possible behaviors of the semi-minimal analysis has been of great interest for pure model theory as well as for its applications, via the $\omega$-stable theory $\mathrm{DCF}_0$, to ordinary algebraic differential equations. See for example \cite{buechler2008vaught}, \cite{jin2018constructing} or \cite{eagles2025splitting}.

Minimal types come equipped with the pregeometry induced by algebraic closure, making their structure simpler to understand than that of arbitrary types (although it can still be very complex). These two results show the crucial importance of minimal types in understanding the structure of finite rank stable theories. It is natural to ask about the connection between the two decompositions, and clarifying this is the main objective of this article. Note that Buechler's levels (see \cite[Section 3]{buechler2008vaught}) are a first answer to this. However, he only considers the specific analysis provided by levels. We work with arbitrary semi-minimal analysis, which should allow our methods to be applied to semi-minimal analyses appearing in nature, such as when studying primitives of solutions of differential equations \cite{jaoui2023relative} or Pfaffian chains \cite{freitag2021not}. Also remark that levels have been generalized to simple theories by Palac\'{i}n and Wagner \cite{palacin2013ample}, and our work here might also generalize.

We will mostly focus on individual steps of the semi-minimal analysis, i.e. fix some type $p \in S(A)$ and some $A$-definable map $f$ with domain containing the realizations of $p$. The image of the realizations of $p$ also is a complete type, which we call $f(p)$. As mentioned previously, we denote such a map $f : p \rightarrow f(p)$. 

A first question one could ask is if the minimal types in the semi-minimal analysis entirely determine those appearing in the domination decomposition, and vice-versa. Neither of these is correct, and examples were already given in \cite{jaoui2023relative}, although not phrased in the language of domination-decomposition. In fact, the notion of \emph{uniform} almost internality from that article is key. It states that the fibers of a map $f$ are all almost internal to a fixed type $r$, and that we can pick the same extra parameters $D$ for all fibers at once. We prove the following, which is Theorem \ref{theo: unif-int-dom-char}:

\begin{customthm}{A}\label{theo: A}
    Let $p \in S(A)$ be stationary and $f : p \rightarrow f(p)$ be an $A$-definable map, and assume that for some $b \models p$, the type $\tp(b/f(b)A)$ is almost internal to some minimal type $r$. Then $p$ is domination equivalent to $f(p) \otimes r^{(m)}$, for some $0 \leq m \leq U(\tp(b/f(b)A))$ (with the convention that if $m = 0$, then $p$ is domination equivalent to $f(p)$). Moreover, if $r \in S(A)$, then $m = U(\tp(b/f(b)A))$ if and only if the map $f$ is uniformly almost $r$-internal.
\end{customthm}

Examples where $m = 0$ are already known from the literature, such as in \cite[Lemma 4.2]{chatzidakis2015differential} and \cite[Theorem 5.2]{jaoui2023relative}, or \cite[Example 6.3.1]{Buechler_2017} in a different context, and show that neither the domination-decomposition nor the semi-minimal analysis determine the other. However, the semi-minimal analysis does determine what types \emph{can} appear in the domination-decomposition, but not with what Morley power, as stated in the following, which is Corollary \ref{cor: from-sma-to-dom-dec} below:

\begin{corollary*}
    Consider a semi-minimal analysis
    \[p=p_n\xrightarrow{f_{n-1}}p_{n-1}\xrightarrow{f_{n-2}}\cdots p_2 \xrightarrow{f_1}p_1 \]
    \noindent and minimal types $r_0, \cdots , r_{n-1}$ such that
    \begin{itemize}
        \item $p_1$ is almost $r_0$-internal and has $U$-rank $n_0$,
        \item for any $1 \leq i \leq n-1$ and any $b_i \models p_i$, the strong type $\stp(b_i/f_{i-1}A)$ is almost $r_i$-internal and has $U$-rank $n_i$.
    \end{itemize}
    Then $p$ is domination equivalent to $ r_0^{(n_0)}\otimes r_1^{(m_1)} \cdots \otimes r_{n-1}^{(m_{n-1})}$ where $0 \leq m_i \leq n_i$ (again with the convention that if $m_i$ is zero then $r_i$ does not appear). 
\end{corollary*}

Given a definable map $f : p \rightarrow f(p)$ with almost $r$-internal fibers, Theorem \ref{theo: A} shows that $p$ is domination equivalent to $f(p) \otimes r^{(m)}$, but does not specify the value of $m$. It is easy to construct examples where any possible value of $m$ is attained. Therefore, isolating conditions in which there is a dichotomy between $m = 0$ and $m$ being maximal, i.e. the $U$-rank of some (any) fiber of $f$, would be useful in concrete domination-decomposition computations. We prove that such a dichotomy is true when the fibers of $f$ have \emph{no proper fibrations} in the sense of \cite{moosa2014some}: a type $q \in S(B)$ has no proper fibration if whenever $a \models q$ and $c \in \dcl(aB)\setminus \acl(B)$, then $a \in \acl(cA)$. 

To show this result, our main tool is Lemma \ref{lem: map-from-no-dom-fix} which allows us, from some type $p \in S(A)$ dominating $r^{(n)}$, for some minimal $r$, to construct a definable map to an $\calR$-internal type, where $\calR$ is the family of $A$-conjugates of $r$. As was pointed out to us by Rahim Moosa, this also yields an interpretation of the \emph{reduction} of the type $p$. Recall that given a family of types $\calQ$ and some type $p$, the $\calQ$-reduction of $p$ is the maximal almost $\calQ$-internal image of $p$ under an $A$-definable map (see \cite{jin2018constructing} for example). 

Let $p \in S(A)$ be a non-algebraic stationary type, domination equivalent to $r_1^{(n_1)} \otimes \cdots \otimes r^{(n_l)}$, where the $r_i$ are pairwise orthogonal minimal types (with $n_i > 0$). We prove that $n_i$ is the $U$-rank of the $\calR_i$-reduction of $p$ over $A$, where $\calR_i$ is the family of $A$-conjugates of $r_i$ over $A$. Note that this can be seen as a type-by-type version of Buechler's work on levels, more precisely \cite[Proposition 3.1]{buechler2008vaught}. However, our proof is independent of his, and our statement more precise.

Note that even though this is a satisfying theoretical result, it does not seem to facilitate the computation of the domination decomposition, mainly because computing the reductions of a type is a difficult problem. For example, there has recently been a flurry of activity around determining which algebraic ordinary differential equations have a minimal generic type in differentially closed fields of characteristic zero. Despite striking recent progress, see for example \cite{freitag2023bounding}, \cite{devilbiss2023generic}, \cite{freitag2023equations}, \cite{jaoui2023density} and \cite{duan2025algebraic}, this remains a challenging problem. 

However, we do remark in Theorem \ref{theo: red-for-autonomous} that for generic types of systems of differential equations of the form:
\[\begin{cases}
    y_1' = f_1(y_1, \cdots , y_k) \\
    \vdots \\
    y_k' = f_k(y_1, \cdots , y_k) 
\end{cases}\]
\noindent where the $f_i$ are rational functions with \emph{constant} coefficients, a slight adaptation of recent work of the authors \cite{eagles2024internality} yields a method to compute their $r$-reduction, where $r$ is the generic type of the constant field.

Given an $A$-definable map $f : p \rightarrow f(p)$ with semi-minimal fibers, another way to obtain a domination-decomposition dichotomy is to make some assumption on the geometry of the fibers. In the last section of this article, we show such a dichotomy in Theorem \ref{theo: trivial-fibers-dicho}, if the fibers are internal to a disintegrated type:

\begin{customthm}{B}\label{theo: B}
    Let $p \in S(A)$ be a stationary type and $f : p \rightarrow f(p)$ be an $A$-definable map such that each fiber is almost internal to a minimal disintegrated type. Then one of the following mutually exclusive possibilities must be true:

    \begin{enumerate}[(a)]
        \item for any $(f(a),f(b)) \models f(p)^{(2)}$, we have $\fibre{p}{a} \perp \fibre{p}{b}$,
        \item there exists a minimal disintegrated type $r \in S(A)$ such that $p$ is interalgebraic with $f(p) \otimes r^{(n)}$, where $n$ is the $U$-rank of any fiber of $f$.
    \end{enumerate}
\end{customthm}

This phenomenon is behind Freitag and Nagloo's recent work \cite{freitag2025algebraic} on relations between solutions of solutions of Painlevé equations of the same family.

Finally, it was noticed while working on this article that the theory is much smoother if one can assume that the types appearing in the domination decomposition of $p \in S(A)$ are themselves over $A$. This is not true in general: in \cite[Proposition 19]{hrushovski1993new} Hrushovski gives an example of an $\aleph_1$-categorical theory with no minimal type over $\acl^{\mathrm{eq}}(\emptyset)$, and in particular any $p \in S(\acl^{\mathrm{eq}}(\emptyset))$ will not have minimal types over $\acl^{\mathrm{eq}}(\emptyset)$ in its domination-decomposition. This phenomenon is crucial in the recent construction by Baldwin, Freitag and Mutchnik \cite{baldwin2024simple} of $\aleph_1$-categorical theories having types with arbitrarily high degree of nonminimality. However, we give some mild conditions on either the ambient theory $T$ or the semi-minimal analysis of $p$ for this to be true. In particular, it will hold in the theories $\mathrm{DCF}_0$ of differentially closed fields of characteristic zero and $\mathrm{CCM}$ of compact complex manifolds. 

We now briefly summarize the contents of this article. Section \ref{sec: prel} contains some preliminaries on internality, domination and orthogonality that we will need. 

Section \ref{sec: dom-and-fib} contains the bulk of the new results of this article. In Subsection \ref{sub: remarks}, we give some basic remarks on the perspective of this work. In particular, we give some conditions guaranteeing that the types appearing in the domination-decomposition of some $p \in S(A)$ are themselves over $A$. In Subsection \ref{sub: fib-appear}, we give results connecting a definable map $f : p \rightarrow f(p)$ to the domination decomposition of $p$, and establish the connection with uniform internality. In Subsection \ref{sub: dom-and-red} we prove the theorem connecting domination decomposition and reductions, and in Subsection \ref{sub: no-prop-fib} we prove the dichotomy result in the case where the fibers of $f$ have no proper fibrations. 

Finally, Section \ref{sec: disintegrated} contains the proof of the splitting dichotomy, when the fibers are almost internal to some disintegrated minimal type.

\section{Preliminaries}\label{sec: prel}

This article is concerned with geometric stability theory, and we will assume that the reader is familiar with it. A good reference is \cite{pillay1996geometric}. We will always work in a superstable theory $T$ in a countable language, which we assume eliminates imaginaries. We will always work inside a fixed, large saturated model $\calU \models T$, and we assume that all parameter sets $A\subset \mathcal{U}$ are small with respect to the saturation of $\mathcal{U}$.

As previously stated, our goal in this article is to explore the relation between two different decompositions of a type in geometric stability: its domination equivalence decomposition and its semiminimal analysis. We will also need various tools such as internality and orthogonality. We recall the relevant definitions and results here. 

We first recall the definition of domination. Note that we follow the convention of Pillay in \cite[Chapter 1, Section 4.3]{pillay1996geometric}: what we call domination for types is also sometimes called eventual domination.

\begin{definition}
    Let $a,b$ be some tuples, and $D$ be some set of parameters. We say that $a$ \emph{dominates $b$ over $D$}, and write $a \trianglerighteq_D b$, if for any set $E$, if $a \forkindep_D E$, then $b \forkindep_D E$. 

    If $p, q \in S(D)$, we say that $p$ dominates $q$ over $D$, and write $p \trianglerighteq_D q$, if there are $a \models p$ and $b \models q$ such that $a \trianglerighteq_D b$.

    If $p,q$ are two stationary types over potentially different sets of parameters, we say that \emph{$p$ dominates $q$}, and write $p \trianglerighteq q$, if there is some $D$ containing the parameters for both $p$ and $q$ such that $p\vert_D \trianglerighteq_D q\vert_D$

    If $p \trianglerighteq q$ and $q \trianglerighteq p$, we will write $p \domequiv q$, and say that $p$ and $q$ are \emph{domination equivalent} (we will also use $\domequiv_D$ in the obvious way).
\end{definition}

Our main reference for domination is \cite{pillay1996geometric}. Domination on types has several nice properties including reflexivity and transitivity. Note that Pillay uses the framework of a-models in that book, which is not necessary for us. Indeed, since we work in a superstable theory $T$ in a countable language, any $\aleph_1$-saturated model is an a-model. We will work with the former for simplicity. Passing to $\aleph_1$-saturated models, one does not need to introduce extra parameters to witness domination. More precisely, we have the following well-known fact, of which we will make frequent use: 

\begin{fact}\label{fact: dom-over-sat}
    Let $p,q$ be any two stationary types. Then $p \trianglerighteq q$ if and only if for some (any) $\aleph_1$-saturated model $M$ containing the parameters of both $p$ and $q$, there are $a \models p\vert_M$ and $b \models q\vert_M$ such that $a \trianglerighteq_M b$.
\end{fact}

For a proof, see for example \cite[Proposition 5.6.4]{Buechler_2017}.

It is well-known that domination is compatible with Morley products (see for example \cite[Remark 5.6.3]{Buechler_2017}:

\begin{fact}\label{fact: dom product} 
Let $p\in S(A)$, $s\in S(C)$ and $r\in S(C')$ be stationary types and suppose that $s\trianglerighteq r.$ Then $p\otimes s\trianglerighteq p\otimes r$.
\end{fact}

One of the most important facts about domination is the well-known \emph{domination-decomposition}, (see for example \cite[Corollary 1.4.5.7]{pillay1996geometric}):

\begin{fact}\label{fact: dom-decomp}
    Let $p \in S(A)$ be stationary of finite Lascar rank. Then for some (any) $\aleph_1$-saturated model $M$ containing $A$, there are pairwise orthogonal minimal types $r_1, \cdots , r_n \in S(M)$ and $k_1, \cdots , k_n \in \mathbb{N}$ such that $p\vert_M \domequiv_M r_1^{(k_1)} \otimes \cdots \otimes r_n^{(k_n)}$. 

    In particular, we have that $p \domequiv r_1^{(k_1)} \otimes \cdots \otimes r_n^{(k_n)}$.
\end{fact}

It is known that the $r_i$ are unique up to ordering and domination-equivalence.

We will use the convention that if some $k_i =0$, then $p \not\trianglerighteq r_i$. Equivalently, the type $r_i$ does not appear in any domination decomposition. 

Note that we stated it here for types of finite rank, but this is not an essential restriction: one could replace minimal types with regular types and obtain the same result.

A useful consequence is the following lemma, stating that minimal types are enough to witness eventual domination.

\begin{lemma}\label{lem: min-dom}

Let $p \in S(A)$ and $q \in S(B)$ be finite $\mathrm{U}$-rank stationary types. Then $p \trianglerighteq q$ if and only if for some (any) $\aleph_1$-saturated model $M \supset A \cup B$, there exist $a \models p\vert_M$ and $b \models q \vert_M$ such that for all $e$ with $\tp(e/M)$ is minimal, we have that $a \forkindep_M e$ implies $b \forkindep_M e$.

\end{lemma}

\begin{proof}
    The left to right direction is immediate. For the right to left direction, suppose that $p \not\trianglerighteq q$. Without loss of generality, assume that $p,q \in S(M)$ by taking non-forking extensions. Fix $a \models p $ and $b \models q$, we need to find $e$ with $\tp(e/M)$ minimal and satisfying $e \forkindep_M a$ and $e \not\forkindep_M b$.

    As $M$ is $\aleph_1$-saturated, we know that there are minimal, pairwise orthogonal types $r_1 , \cdots , r_n \in S(M)$ and $k_1 , \cdots , k_n \in \mathbb{N}$ such that $q \domequiv_M r_1^{(k_1)} \otimes \cdots \otimes r_n^{(k_n)}$. For all $i$, let $l_i := \max \{j \in \mathbb{N}: p \trianglerighteq_M r_i^{(j)} \}$ (note that $l_i$ could be $0$). There must be $i$ such that $l_i < k_i$. Else, as $p \trianglerighteq_M r_i^{(l_i)}$ for all $i$, and the $r_i$ are pairwise orthogonal, we get $p \trianglerighteq_M r_1^{(l_1)} \otimes \cdots \otimes r_n^{(l_n)}$. As $l_i \geq k_i$ for all $i$, this would imply $p \trianglerighteq_M q$.

    Without loss of generality, assume that $l_1 < k_1$, let $k = k_1$, $l=l_1$ and $r_1 = r$. Let $(c_1, \cdots, c_k) \models r^{(k)}$ be such that $b \trianglerighteq_M c_1 \cdots c_k$. In particular we see that $c_i \not\forkindep_M b$ for all $i \in \{ 1 , \cdots , k\}$. 
    There must be $i \in \{ 1 , \cdots , k \}$ such that $a \forkindep_M c_i$. Else by minimality $c_i \in \acl(aM)$ for all $i$, and thus for any tuple $d$, if $d \forkindep_M a$, then $d \forkindep_M c_1 \cdots c_k$, which implies $a \trianglerighteq_M c_1 \cdots c_k$ and thus $p \trianglerighteq_M r^{(k)}$, contradicting $l<k$. Thus $e = c_i$ is the element we were looking for.
\end{proof}

\begin{remark}
    From the proof, we see that we could pick $e$ to be a realization of $r_i$, for any $r_i$ appearing in the domination decomposition of $p$ with a higher Morley power than in the decomposition of $q$. 
\end{remark}

Let us now recall the notions of internality and analysis, as well as some basic results.

\begin{definition}
    Let $p \in S(A)$ be a stationary type and $r \in S(B)$. We say that $p$ is (resp. almost) $r$-internal if there is $D \supset A \cup B$, some $a \models p \vert_D$ and some $c_1, \cdots , c_n \models r$ such that $a \in \dcl(c_1, \cdots , c_n D)$ (resp. $\acl$).
\end{definition}

A type that is almost internal to a minimal type is called \emph{semiminimal}.

We will also consider internality to invariant families of types. Let $\calR$ be a a family of partial types, potentially over different parameters. A \emph{realization} of $\calR$ is a tuple realizing any of the partial types in $\calR$. We say it is a realization of $\calR$ over $B$ if it is a realization of a partial type in $\calR$ with parameters in $B$. We say that $\calR$ is $A$-\emph{invariant} if for any $\sigma \in \mathrm{Aut}(\calU)$ fixing $A$ pointwise and any realization $c$ of $\calR$, the tuple $\sigma(c)$ also is a realization of $\calR$. 

The two kinds of invariant families we will consider in this document are $\{ r \}$ for some type $r \in S(A)$, and, given some type $s$ over any parameters, the family $\calS$ of its conjugates under the group $\mathrm{Aut}_A(\calU)$ of automorphisms of $\calU$ fixing $A$ pointwise. 

Fix an $A$-invariant family of partial types $\calR$. We recall the definition of (almost) internality in that context:

\begin{definition}
    Let $p \in S(A)$ be stationary, we say $p$ is $\calR$-internal (resp. almost $\calR$-internal) if there are $a \models p$, some $B \supset A$ with $a \forkindep_A B$ and some realizations $c_1, \cdots , c_n$ of $\calR$ over $B$ such that $a \in \dcl(c_1, \cdots , c_n , B)$ (resp. $\acl$).
\end{definition}

When $\calR = \{ r \}$ is given by a single type over $A$, we recover (almost) $r$-internality. 

The presence of the extra parameters $B$ in the definition is crucial. However, it is well-known that in the case of internality to a disintegrated minimal type, almost internality collapses to algebraicity. Since we could not locate a proof in the literature, we provide one. Recall, a minimal type $r\in S(A)$ is \textit{disintegrated} if for every subset $X\subseteq r(\calU),$ $\acl(X)=\bigcup\limits_{x\in X}\acl(x).$ Sometimes these types are referred to as trivial minimal types. 

\begin{lemma}\label{lem: int-to-disin-implies-alg}
    Let $B \supset A$ and $p \in S(B)$ be stationary and almost internal to some minimal, disintegrated type $r \in S(A)$. Then $p$ is $r$-algebraic, meaning that for any $a \models p$, there are $c_1, \cdots, c_n \models r$ such that $a \in \acl(c_1, \cdots , c_n , B)$.
\end{lemma}

\begin{proof}
    Notice that in fact $p$ is $r \vert_B$ almost internal. Indeed, if $a \models p$ is such that $a \in \acl(c_1, \cdots , c_n , e , B)$ for some $c_i \models r$ and parameters $e$, then any $c_i$ not independent from $B$ over $A$ is in $\acl(B)$, and thus can be removed. Therefore, replacing $r$ with $r\vert_B$, we may assume that $B = A$.

    There is a type $q \in S(B)$ such that $q$ is $r$-internal and $p$ and $q$ are interalgebraic (see for example \cite[Lemma 3.6]{jin2020internality}), meaning that for any $a \models p$, there is $c \models q$ such that $\acl(aB) = \acl(cB)$. Consider the binding group $\Aut_B(q/r(\mathcal{U}))$, which is the group of permutations of $q(\mathcal{U})$ induced by automorphisms of $\mathcal{U}$ fixing $B \cup r(\mathcal{U})$ pointwise. It is well known that this group is type-definable, and definably isomorphic to a group type-definable in $r^{\mathrm{eq}}$ (see, for example, Theorem 7.4.8 and Remark 7.4.9 of \cite{pillay1996geometric}).

    As $r$ is disintegrated, there cannot be any infinite group type-definable in $r^{\mathrm{eq}}$. This is well-known, but here is a quick justification: disintegratedness implies that $r$, and thus also $r^{\mathrm{eq}}$, has \emph{trivial forking}: any three pairwise independent $a,b,c$ in $r^{\mathrm{eq}}$ must in fact be independent. However, if some infinite group $(G, \cdot)$ is type-definable, picking independent generics $a,b$, then $a,b$ and $a \cdot b$ are pairwise independent, but not independent.
    
    In particular $\Aut_B(q/r(\mathcal{U}))$ is finite. Hence $q$, and therefore $p$, are $r$-algebraic.
\end{proof}

Semiminimal types give rise to another decomposition for finite Lascar rank types in a superstable theory, called the \emph{semiminimal analysis}. Instead of Morley products, it uses \emph{fibrations}:

\begin{definition}
    Let $p \in S(A)$. A definable map on $p$ is an $A$-definable partial map $f$, with domain containing $p$. Note that the image $f(p(\calU))$ is a complete type, which we denote $f(p)$. We denote such a definable map by $f : p \rightarrow f(p)$.
    
    If for any $a \models p$, the type $\tp(a/f(a)A)$ is stationary, we say that $f$ is a fibration. We will often denote this type by $\fibre{p}{a}$, and call it the \emph{fiber} of $f$.     
\end{definition}

Note that we adopt a slightly different point of view as the one from the introduction, which allows us to treat all steps of the analysis uniformly.

\begin{fact}
    Let $p \in S(A)$ be a finite Lascar rank stationary type, and pick any $a \models p$. 
    
    There are $a_0,\cdots,a_n$ with $a_n:=a$, $a_0=\emptyset$ such that $a_i\in\dcl(a_{i+1}A)$ for all $i=1,\cdots,n-1$, and $\stp(a_{i+1}/a_iA)$ is semiminimal for all $i=0,\cdots,n-1$. 
    
    This gives rise to a sequence of types and $A$-definable functions \[p=p_n\xrightarrow{f_{n-1}}p_{n-1}\xrightarrow{f_{n-2}}\cdots\xrightarrow{f_1}p_1\xrightarrow{f_0}\bullet\] where $p_i=\tp(a_i/A)$ and $\bullet$ represents some element in $A$ (so $f_0$ has only one fiber). We call this sequence of types and functions a semiminimal analysis of $p$ and $n$ the length of this semiminimal analysis.
\end{fact}

One of the main questions we ask in this article is the connection between the semi-minimal analysis and the domination-decomposition. It is often enough to look at one step of the semi-minimal analysis. 

We end this section with some well-known results connecting internality, orthogonality and domination equivalence.

This is \cite[Chapter XIII, Corollary 2.24]{baldwin2017fundamentals}:

\begin{fact}\label{fact: min ortho trans}
    Let $p \in S(A)$, $q \in S(B)$ and $r \in S(C)$ be stationary types, with $r$ minimal. If $p \not\perp r$ and $r \not\perp q$, then $p \not\perp q$.
\end{fact}

The following well-known connection between internality and domination will be generalized in this article:

\begin{proposition}\label{pro: domeq for internal}
    Let $p \in S(A)$ be stationary non-algebraic and let $r \in S(A')$ be a minimal type. Also fix $n := \mathrm{U}(p) > 0$. The following are equivalent:
    \begin{enumerate}
        \item $p$ is almost $r$-internal,
        \item there exist some set of parameters $B \supset A\cup A'$, some $a \models p\vert_B$ and some $c_1, \cdots , c_n \models r^{(n)}\vert_B$, such that $a$ is interalgebraic, over $B$, with $c_1 , \cdots , c_n$,
        \item $p \domequiv r^{(n)}$.
    \end{enumerate}
    
\end{proposition}
    \begin{proof}
    
    We first prove that $(1) \Rightarrow (2)$. So assume that $p$ is almost $r$-internal. The extension $p|_{A'}$ is still almost $r$-internal. The assumption gives us some $B\supset A'\cup A$, some $a\models (p|_{A'})|_B$, and $c_1, \cdots ,c_m$ realizations of $r$, such that $a\in\acl(c_1, \cdots ,c_m,B)$. We can pick $B$ such that $m$ is minimal (among all possible such $m$ and $B$). As the type $r$ is minimal, this implies that $(c_1, \cdots ,c_m)\models (r|_B)^{(m)}$. 

    It is enough to show that $c_i\not\forkindep_Ba$ for all $i$. If on the contrary $c_1\forkindep_Ba$ (without loss of generality), as $c_1 \cdots c_m \not\forkindep_B a$, we get $c_2 \cdots c_m\not\forkindep_{Bc_1}a.$ As $c_1\forkindep_Ba$ and $a\forkindep_{A'}B,$ we get $a\forkindep_{A'}c_1B$, which contradicts the minimality of $m$.

    So $a$ is interalgebraic over $B$ with $(c_1, \cdots , c_m)\models(r^{(m)})|_B$. This implies $m=n$.

    The implication $(2) \Rightarrow (3)$ is immediate, so we now prove $(3) \Rightarrow (1)$. Assume that $p \domequiv r^{(n)}$. Then there are $B \supset A \cup A'$, some $a \models p\vert_B$ and some $(c_1, \cdots , c_n) \models r^{(n)}\vert_B$ such that $a \domequiv_B (c_1, \cdots , c_n)$. We will prove that $a \not\forkindep_{c_1 \cdots c_{i-1}} c_i$ for all $1 \leq i \leq n$, which implies $a \in \acl(c_1 \cdots c_n B)$ because $n = U(p)$, and since $a \forkindep_A B$, that $p$ is almost $r$-internal. Note first that since $c_i \not\forkindep_B c_1 \cdots c_n$ for all $i$, we have $a \not\forkindep_B c_i$ for all $i$. Assume, for a contradiction, that $a \forkindep_{c_1 \cdots c_{i-1}B} c_i$. Then as $c_i \forkindep_B c_1 \cdots c_{i-1}$, we obtain $a \forkindep_B c_i$, a contradiction. Therefore $a \not\forkindep_{Bc_1 \cdots c_{i-1}} c_i$ for all $i$, so $a \in \acl(c_1 \cdots c_n B)$.
\end{proof}

As a corollary, we obtain:

\begin{corollary}\label{cor: min ortho int}

Let $p \in S(A),$ $r\in S(A')$ and $q \in S(B)$ be stationary types with $A \subset B$. Suppose that $r$ is minimal, $q$ is non-algebraic, and $q$ is almost $r$-internal. If $p \perp r$ then $p \perp q$.

\end{corollary}

\begin{proof}
    By Proposition \ref{pro: domeq for internal}, we have that $q \domequiv r^{(n)}$, where $n > 0$ is the $U$-rank of $q$. Assume that $p \perp r$. Then we also have $p \perp r^{(n)}$, so $p \perp q$. 
\end{proof}

\section{Domination and fibrations}\label{sec: dom-and-fib}

In this section, we explore the connection between fibrations and domination decomposition. We start by, in the next subsection, justifying some of the assumptions we will make on our fibrations.

\subsection{Remarks on the assumptions}\label{sub: remarks}

Recall that one of our stated goals is to understand the connection between the domination equivalence decomposition of a type and its semiminimal analysis. We will often focus on one step of the semiminimal analysis, i.e. we will consider a type $p \in S(A)$ and some $A$-definable fibration $f : p \rightarrow f(p)$ with semiminimal fibers. By the fibers of a fibration, we mean the stationary types $\fibre{p}{b} = \tp(b/f(b)A)$, for some $b \models p$. 

The main question we want to answer is whether or not the fibers of $f$ appear in the domination decomposition of $p$. More precisely, we always have that $p \trianglerighteq f(p)$, and by the fibers appearing, we mean that $f(p) \not\trianglerighteq p$, in other words, that $p$ and $f(p)$ are not domination equivalent. 

From that perspective, we give a sufficient condition for the fibers to not appear:

\begin{proposition}
    Let $p\in S(A)$ be stationary of finite $U$-rank and $f:p\to f(p)$ a fibration with semiminimal fibers. If $p\perp \fibre{p}{b}$ for some (any) $f(b)\models f(p)$, then $p\domequiv f(p)$.
    \label{prop: ortho-to-fib}
\end{proposition}
\begin{proof}
    Let $M$ be an $\aleph_1$-saturated model containing $A$ and $b\models p \vert _M.$ Assume that $f(p) \not\trianglerighteq p$, which implies in particular $f(b) \not\trianglerighteq_M b$.
    
    By Lemma \ref{lem: min-dom}, there is some $e$ with $r = \tp(e/M)$ minimal such that $f(b) \forkindep_M e$ and $b \not\forkindep_M e$. This implies that $b \not\forkindep_{f(b) M} e$, so $\tp(b/f(b)M) \not\perp r$. Note that $\tp(b/Mf(b))=(p|_M)_{f(b)}=(\fibre{p}{b})|_M$, thus $\fibre{p}{b}\not\perp r$. As $r\not\perp p$, transitivity of non-orthogonality for minimal types (i.e. Fact \ref{fact: min ortho trans}) gives $\fibre{p}{b}\not\perp p$.
\end{proof}

Therefore, we will often assume that $p \not\perp \fibre{p}{b}$, as the other case is completely solved by the previous lemma. This has the consequence that the fibers are pairwise nonorthogonal:

\begin{proposition}\label{pro: non-orth-to-fib}
    Let $p \in S(A)$ be stationary, and some fibration $f : p \rightarrow f(p)$. If $p \not\perp \fibre{p}{b}$ for some (any) $f(b) \models f(p)$, then for any independent $f(b), f(c) \models f(p)$, we have $\fibre{p}{b} \not\perp \fibre{p}{c}$. 
    
    Moreover, if each fiber $\fibre{p}{b}$ is almost internal to some minimal type $r_{f(b)}$, then for any $f(b),f(c) \models f(p)$, we have both $r_{f(b)} \not\perp r_{f(c)}$ and $\fibre{p}{b} \not\perp \fibre{p}{c}$.
\end{proposition}

\begin{proof}
    Let $f(b) \models f(p)$ be such that $p \not\perp \fibre{p}{b}$. Then $\fibre{p}{b}$ is \emph{non orthogonal to $A$} in the sense of \cite[Lemma 4.3.3]{pillay1996geometric}: some type $q \in S(B)$, for some set of parameters $B$, is non-orthogonal to $A$ is it is non-orthogonal to some type over $\acl(A)$. By that same lemma, for any $f(b),f(c) \models f(p)^{(2)}$, we have that $\fibre{p}{b} \not\perp \fibre{p}{c}$. 
    
    Now we prove the moreover part, assume each fiber $\fibre{p}{b}$ internal to some minimal type $r_{f(b)}$, over some parameters $A_{f(b)}$ containing $f(b)$. Let $f(b), f(c)$ be any two realizations of $f(p)$. Pick some $f(d) \models f(p)\vert_{A_{f(b)},A_{f(c)}}$. By the first part of the lemma, we have $\fibre{p}{b} \not\perp \fibre{p}{d}$ and $\fibre{p}{d} \not\perp \fibre{p}{c}$, and thus $r_{f(b)} \not\perp r_{f(d)}$ and $r_{f(d)} \not\perp r_{f(c)}$ by Corollary \ref{cor: min ortho int}. This gives $r_{f(b)} \not\perp r_{f(c)}$ by Fact \ref{fact: min ortho trans}, which finally implies that $\fibre{p}{b} \not\perp \fibre{p}{c}$ by Fact \ref{fact: min ortho trans}.
\end{proof}

The converse of Proposition \ref{pro: non-orth-to-fib} is false. See Example \ref{ex: JJP} for a fibration $f : p \rightarrow f(p)$ such that for any $f(b),f(c)$ we have $\fibre{p}{b} \not\perp \fibre{p}{c}$ but $p \perp \fibre{p}{b}$. 

What the proposition tells us is that at each step of the analysis where the fibers appear in the domination decomposition, we obtain a family of pairwise non-orthogonal minimal types, potentially over extra parameters. There is in general no canonical choice for which to pick for the domination decomposition. This has the effect of making many statements very cumbersome. We would like to be able to assume that the fibers are all non-orthogonal to some minimal type over $A$. In the rest of this section, we point out two conditions under which this is true. 

The first is if the fibers are internal to locally modular types. More generally:

\begin{lemma}\label{lem: min-descent-loc-mod}
    Let $p = \tp(a/A) \in S(A)$ be stationary, and suppose that $p \not\perp t$, for some minimal, locally modular $t$, potentially over some extra parameters. Then there exists a minimal $r \in S(A)$ such that $p \not\perp r$ and $t \domequiv r$.
\end{lemma}

\begin{proof}
    This is the same as the proof of \cite[Proposition 2.3]{moosa2014some}, but we give the full proof for the comfort of the reader. Let $\mathcal{T}$ be the set of $A$-conjugates of $t$. This set is $A$-invariant, and $p$ is not foreign to $\mathcal{T}$ (in the sense of \cite[Chapter 7, 4.1]{pillay1996geometric}). Therefore, by \cite[7.4.6]{pillay1996geometric}, there is $c \in \dcl(aA) \setminus \acl(A)$ such that $\stp(c/A)$ is $\mathcal{T}$-internal. Since $c \in \dcl(aA)$ and $\tp(a/A)$ is stationary, the type $q = \tp(c/A)$ is also stationary, and in particular $\tp(c/A) \models \tp(c/\acl(A))$. Therefore $q$ is also $\mathcal{T}$-internal. This in fact implies that $q$ is almost $t$-internal (see \cite[Lemma 5.4]{freitag2025finite} for a proof of that well-known, but non-trivial, fact).
    As $t$ is locally modular, the type $q$ is one-based (this is not immediate, it was proven by Wagner in \cite[Corollary 9]{wagner2004some}). If it has U-rank one, then it is minimal and we are done. Else, consider some tuple $e$ such that $U(c/Ae) = U(c/A) - 1$ and $e = \Cb(c/Ae)$. By one-basedness, we have $e \in \acl(Ac)$, and thus:
    \[U(c/A) = U(e/Ac) + U(c/A) = U(ce/A) = U(c/eA) + U(e/A) = U(c/A) + U(e/A) - 1\]
    \noindent so $U(e/A) = 1$. Moreover, we have that $e \in \dcl(c_1, \cdots , c_l)$, for $c_1, \cdots ,c_l$ some Morley sequence in $\stp(c/Ae)$, which implies that the type $r = \tp(e/A)$ is stationary, and thus minimal. As $e \in \acl(cA) \setminus \acl(A)$ and $c \in \dcl(aA)$, we get $r \not\perp q$ and $r \not\perp p$. As $q$ is almost $t$-internal, Corollary \ref{cor: min ortho int} gives us that $r$ is nonorthogonal, and thus domination equivalent, to $t$. 
\end{proof}

Given this lemma, we can obtain what we want under the following:

\begin{Assumption}\label{assump: non-orth}

There is a non-locally modular minimal type $ \mathfrak{f} \in S(\emptyset)$ such that any other minimal non-locally modular type is non-orthogonal to $\mathfrak{f}$.

\end{Assumption}

This is a frequent assumption in the geometric stability literature (see for example \cite{palacin2017definable}) and is satisfied by many theories of interest: it is true in $\mathrm{DCF}_{0,m}$ and $\mathrm{CCM}$. Note that some more unusual theories do satisfy this properties, for example the counterexample to the canonical base property of Hrushovski, Palac\'{i}n and Pillay in \cite{hrushovski2013canonical} (see also \cite{loesch2022possibly} for more examples). The type $\mathfrak{f}$ does not even have to be the generic type of a field, as illustrated by the counterexamples to the canonical base property based on Baudisch's group \cite{baudisch1996new}, see \cite{blossier2022cm} and \cite{loesch2025rings}.

For an example of a theory not satisfying the assumption, we also follow \cite{palacin2017definable}: consider a two sorted theory with one sort for a pure set $I$, and another sort $S$ with a surjection $\pi : S \rightarrow I$ such that each fiber is an algebraically closed field of fixed characteristic (and there is no relation between the fibers). 

\begin{lemma}[Under assumption \ref{assump: non-orth}]\label{lem: min descent}

Let $p = \tp(a/A) \in S(A)$ be stationary, and suppose that $p \not\perp t$, for some minimal $t$, potentially over some extra parameters. Then there exists a minimal $r \in S(A)$ such that $p \not\perp r$ and $t \domequiv r$.

\end{lemma}

\begin{proof}

First notice that if $t$ is non-locally modular, by assumption \ref{assump: non-orth}, it must be non-orthogonal, and hence domination equivalent, to $\mathfrak{f}$. If $t$ is locally modular, we can conclude by Lemma \ref{lem: min-descent-loc-mod}.
\end{proof}

Note that the previous example of a theory not satisfying assumption \ref{assump: non-orth}, the conclusion of Lemma \ref{lem: min descent} should still hold. Indeed, the natural candidate to not satisfy \ref{lem: min descent} would be the generic type of the sort $F$, but it is non-orthogonal to the generic type of $I$. 

However, there are theories which do not satisfy the conclusion of Lemma \ref{lem: min descent}, the first of which was given by Hrushovski in \cite[Proposition 19]{hrushovski1993new}, using his groundbreaking constructions. We give a brief explanation. The theory in question is almost strongly minimal of Morley rank 2. Recall that almost strong minimality means that given a saturated model $\calU$, there is a formula $\phi$ (maybe using parameters) such that $\phi(\calU)$ is strongly minimal, and some finite set $A$ with $\calU \subset \acl(A, \phi(\calU))$. Note that this implies in particular that $U$-rank and Morley rank coincide. Hrushovski shows that in this theory, there are no strongly minimal, and hence minimal, type over $\acl^{\mathrm{eq}}(\emptyset)$. Taking any type $p \in S(\acl^{\mathrm{eq}}(\emptyset))$, it thus cannot have any minimal type over $\acl^{\mathrm{eq}}(\emptyset)$ in its domination-decomposition. Note that this property is not without interesting consequences. For example, it plays a crucial role in the recent proof by Baldwin, Freitag and Mutchnik \cite{baldwin2024simple} that the degree of nonminimality is unbounded in $\aleph_1$-categorical theories. 

Back to our main interest, this is what we have obtained in the context of Proposition \ref{pro: non-orth-to-fib}:

\begin{corollary}\label{cor: suff-for-fib-descent}
    Let $p \in S(A)$ be stationary, and some fibration $f : p \rightarrow f(p)$, with each fiber $\fibre{p}{b}$ almost internal to some minimal type $r_{f(b)}$. If $p \not\perp \fibre{p}{b}$ for some (any) $f(b) \models f(p)$ and either:
    \begin{itemize}
        \item some (any) $r_{f(b)}$ is locally modular,
        \item $T$ satisfies Assumption \ref{assump: non-orth},
    \end{itemize}
    \noindent then there exists a minimal $r \in S(A)$ such that any $\fibre{p}{b}$ is almost $r$-internal.
\end{corollary}

\begin{proof}
    Let $f(b)\models f(p)$ be such that $p\not\perp \fibre{p}{b}$. Since $\fibre{p}{b}$ is non-algebraic and almost internal to the minimal type $r_{f(b)}$, we have $\fibre{p}{b}\not\perp r_{f(b)}$, which gives $p \not \perp r_{f(b)}$ by Corollary \ref{cor: min ortho int}. Under either assumption, there is some $r\in S(A)$ such that $\fibre{p}{b}\not\perp r$ and $r\domequiv r_{f(b)}$. Since $\fibre{p}{b}$ is almost internal to $r_{f(b)}$, Proposition \ref{pro: domeq for internal} gives us $\fibre{p}{b} \domequiv r_{f(b)}^{(n)}$, where $n$ is the $U$-rank of $\fibre{p}{b}$. So $\fibre{p}{b} \domequiv r^{(n)}$, hence $\fibre{p}{b}$ is almost $r$-internal, again by Proposition \ref{pro: domeq for internal}. 
    
\end{proof}

Finally, we also have an straightforward consequence for the domination decomposition:

\begin{theorem}[Under assumption \ref{assump: non-orth}]\label{theo: min-descent}

Let $p \in S(A)$ be any finite Lascar rank stationary type. Then there are $r_1, \cdots, r_n \in S(A)$ minimal, pairwise orthogonal types such that $p \domequiv r_1^{(k_1)} \otimes \cdots \otimes r_n^{(k_n)}$ for some $k_1, \cdots ,k_n \in \mathbb{N}$.

\end{theorem}

\begin{proof}
   By Fact \ref{fact: dom-decomp}, there is a model $M \models T$, some $t_1, \cdots, t_n \in S(M)$ minimal types and some $k_1, \cdots , k_n \in \mathbb{N}$ such that $p \domequiv t_1^{(k_1)} \otimes \cdots \otimes t_n^{(k_n)}$, and in particular $p \not\perp t_i$ for all $i$. By Lemma \ref{lem: min descent}, there are minimal types $r_1 , \cdots , r_n \in S(A)$ such that $r_i \domequiv t_i$ for all $i$, and thus $r_i^{(k_i)} \domequiv t_i^{(k_i)}$. We can now get the result by repeatedly applying Fact \ref{fact: dom product} (and noticing that for any types $r,s$, we have $r \otimes s \domequiv s \otimes r$).
\end{proof}
\subsection{When the fibers appear}\label{sub: fib-appear}

We are now going to give more precise results on whether fibers coming from a fibration will appear in the domination decomposition. The main lemma is the following:

\begin{lemma}
\label{lem: dom for internal fiber}
    Let $p\in S(A)$ be stationary. Let $f:p\to q$ be a fibration such that for some (any) $a\models p$, the type $\stp(a/f(a)A)$ is almost internal to a minimal type $r_{f(a)}$. Then for some $0\leq m \leq\mathrm{U}(\stp(a/f(a)A))=n$,  we have $f(p)\otimes r_{f(a)}^{(m)}\domequiv p$.
\end{lemma}

\begin{proof}
    For any $a \models p$, we let $A_{f(a)}$ be such that $r_{f(a)} \in S(A_{f(a)})$. We may assume $\acl(f(a)  A )\subset A_{f(a)}$. 
    
    We fix some $f(a) \models f(p)$ and some $\aleph_1$-saturated model $M$ containing $A_{f(a)}$. If $f(p) \domequiv p$, we let $m = 0$ and we are done. We can thus assume that $f(p) \not\trianglerighteq p$. Fix some $b \models p\vert_M$.
    
    We will inductively construct some $(c_1, \cdots, c_k) \models (r_{f(a)}\vert_M)^{(k)}$ with the following properties:
    \begin{enumerate}
        \item $c_1, \cdots , c_k \in \acl(bM)$,
        \item $c_1 \cdots c_k \forkindep_M f(b)$
        \item $b \not\forkindep_{c_1 \cdots c_{k-1}f(b)M} c_k$.
    \end{enumerate}

    More precisely, we show that, given such $c_1, \cdots , c_k$, as long as $f(p) \otimes r_{f(a)}^{(k)} \not \trianglerighteq p$, we can obtain $c_{k+1}$ such that $c_1, \cdots, c_k , c_{k+1}$ have the same properties. 

    So assume that $f(p) \otimes r_{f(a)}^{(k)} \not \trianglerighteq p$, and fix some $c_1, \cdots , c_k$ with the given properties (this will also work if $k = 0$). By Lemma \ref{lem: min-dom} there is some $e$ such that $r = \tp(e/M)$ is minimal and:
    \begin{enumerate}[(a)]
        \item $e \forkindep_M c_1 \cdots c_k f(b) $,
        \item $e \not\forkindep_M b$.
    \end{enumerate}
    By (a), we see that $e \forkindep_M f(b)$, and therefore (b) gives us $e \not\forkindep_{f(b)M} b$. This implies that $r \not\perp \stp(b/f(b)A)$, and by Corollary \ref{cor: min ortho int} that $r \not\perp r_{f(b)}$. Finally, Propositions \ref{prop: ortho-to-fib} and \ref{pro: non-orth-to-fib} together imply that $r_{f(b)} \not\perp r_{f(a)}$, therefore $r \not\perp r_{f(a)}$. 

    Since $M$ is $\aleph_1$-saturated, this implies that $r \not\perp^w r_{f(a)}\vert_M$, and because the types are minimal, there is some $c_{k+1} \models r_{f(a)}\vert_M$ that is interalgebraic with $e$ over $M$. Conditions (a) and (b) are up to interalgebraicity over $M$, so we may assume that $e = c_{k+1}$, i.e. we have:
    \begin{enumerate}[(a)]
        \item $c_{k+1} \forkindep_M c_1 \cdots c_k f(b) $
        \item $c_{k+1} \not\forkindep_M b$.
    \end{enumerate}
    This immediately yields $c_{k+1} \in \acl(bM)$, so we have condition (1) for $k+1$. Condition (2) for $k$ and (a) give condition (2) for $k+1$. Finally, suppose that condition (3) for $k+1$ was false, so we have $b \forkindep_{c_1 \cdots c_{k}f(b)M} c_{k+1}$. Using condition (a), we obtain that $b \forkindep_{M} c_{k+1}$, a contradiction. 

    At some step $k \leq n$, we must have that $p \domequiv (r_{f(a)})^{(k)} \otimes f(p)$. Indeed, by condition (3), we have $U(\tp(b/c_1 \cdots c_{k}f(b)M)) < U(\tp(b/c_1 \cdots c_{k-1}f(b)M))$ for all $k \geq 1$. 
\end{proof}

Remark that if either $r_{f(a)}$ is locally modular or Assumption \ref{assump: non-orth} is true, we can replace $r_{f(a)}$ by some minimal type $r \in S(A)$.

Given this result, we can precisely spell out the connection between the domination decomposition and semi-minimal analysis. 

\begin{corollary}\label{cor: from-sma-to-dom-dec}
    Let $p \in S(A)$ and
    \[p=p_n\xrightarrow{f_{n-1}}p_{n-1}\xrightarrow{f_{n-2}}\cdots\xrightarrow{f_1}p_1\xrightarrow{f_0}\bullet\]
    \noindent be a semi-minimal analysis. Then $p\domequiv r_0^{(n_0)} \otimes r_1^{(m_1)} \otimes \cdots \otimes r_{n-1}^{(m_{n-1})}$ where $r_i$ is a minimal type such that for some $b_i \models p_i$, the strong type $\stp(b_i/f_{i-1}(b_i)A)$ is almost $r_i$-internal and $0 \leq m_i \leq n_i$ where $n_i$ is the $U$-rank of some (any) fiber for all $i=0,\cdots n-1$. 
\end{corollary}

\begin{proof}
    Fix $a\models p$ and set $a_n:=a.$ Inductively set $a_i=f_i(a_{i+1})$ for all $i=1,\cdots,n-1$ and $a_0=\emptyset.$ We induct on $n$, the length of the analysis. If $n=1$, then $p$ itself is semiminimal and by Proposition \ref{pro: domeq for internal}, $p\domequiv r^{(m)}$ with $m=U(p).$

    Suppose the result holds for all analyses of length $n-1.$ In particular, by inductive hypothesis we have that $p_{n-1}\domequiv r_0^{(m_0)}\otimes \cdots\otimes r_{n-1}^{(m_{n-1})}$ satisfying the statement of the corollary. Let $r_{n-1}$ a minimal type such that $\stp(a_n/a_{n-1}A)$ is almost $r_{n-1}$-internal. By Lemma \ref{lem: dom for internal fiber} we obtain $p\domequiv p_{n-1}\otimes r_{n-1}^{(m_{n-1})}$ for some $0\leq m_{n-1} \leq \mathrm{U}(a_n/a_{n-1}A)$. This gives the result by using Fact \ref{fact: dom product}.
\end{proof}

This corollary makes no extra assumption on either the theory $T$ or the pregeometry of the $r_i$, and therefore the $r_i$ may be over extra parameters. In the rest of this subsection, we will make the assumption that fibers of our fibration are semiminimal, and internal to a minimal type over the base parameters. As explained in Subsection \ref{sub: remarks}, this is always the case if this minimal type is locally modular, and in general under Assumption \ref{assump: non-orth}, which holds in many theories of interest, such as $\mathrm{DCF}_0$ and $\mathrm{CCM}$.

Consider a fibration $f : p \rightarrow f(p)$, and assume that the fibers are almost $r$-internal, for some $r \in S(A)$. We know from Lemma \ref{lem: dom for internal fiber} that $p \domequiv r^{(m)} \otimes f(p)$, for some $0 \leq m \leq n$, where $n$ is the $U$-rank of any fiber. What are the possibilities for $m$? We note that $m = 0 $ is possible, by Proposition \ref{prop: ortho-to-fib}:

\begin{example}[\cite{jaoui2023relative}, Example 5.6]\label{ex: JJP}
    In $\mathrm{DCF}_0$, consider the generic type $p$ of the system:
    \[
    \begin{cases}
        x' = x^3 (x-1) \\
        y ' = xy
    \end{cases}
    \]
    and the fibration $f$ given by projection on the $x$ coordinate. It is shown in \cite[Corollary 5.5]{jaoui2023relative} that $p$ is orthogonal to the field of constants $\calC$. But for any $b \models p$, the fiber $\fibre{p}{b}$ is $\calC$-internal, and therefore $p \perp \fibre{p}{b}$, which implies, by Proposition \ref{prop: ortho-to-fib}, that $p \domequiv f(p)$.
\end{example}

The next result points out exactly when $m$ is maximal: it corresponds to \emph{uniform almost internality} of $f$, a notion defined in \cite{jaoui2023relative}: 

\begin{definition}
    Let $p \in S(A)$ and $f : p \rightarrow f(p)$ be a fibration, and let $r \in S(A)$. We say $f$ is uniformly $r$-internal if there are $b \models p$ and some $D \supset A$ such that:
    \begin{itemize}
        \item $b \in \dcl(c_1, \cdots , c_m, f(b), D)$ for some realizations $c_1, \cdots , c_m$ of $r$,
        \item $b \forkindep_A D$.
    \end{itemize}
    If $b \in \acl(c_1, \cdots , c_m, f(b), D)$ instead, we say that $f$ is uniformly almost $r$-internal.
\end{definition}

We have:

\begin{theorem}\label{theo: unif-int-dom-char}
    Let $p \in S(A)$ be stationary and $f : p \rightarrow f(p)$ be a fibration. Suppose that there is some minimal $r$ (over any small set of parameters) such that for some (any) $b \models p$, the fiber $\fibre{p}{b}$ is almost $r$-internal. Let $n = U(\fibre{p}{b})$ for some (any) $b \models p$. 

     Then $p \domequiv r^{(m)} \otimes f(p)$ for some $0 \leq m \leq n$, and if $r \in S(A)$ the following are equivalent:

     \begin{enumerate}
         \item $f : p \rightarrow f(p)$ is uniformly almost $r$-internal,
         \item there is some $D \supset A$ such that $p \vert_D$ is interalgebraic with $(f(p) \otimes r^{(n)})\vert_D$,
         \item $p \domequiv f(p) \otimes r^{(n)}$.
     \end{enumerate}
     Moreover if $f(p)$ is almost $r$-internal, then these conditions are equivalent to:
     \begin{enumerate}
         \item[(4)] $p$ is almost $r$-internal.
     \end{enumerate}
\end{theorem}

\begin{proof}
    The first part is simply Lemma \ref{lem: dom for internal fiber}. We start by proving $(1) \Rightarrow (2)$. Assume that $f : p \rightarrow f(p)$ is uniformly almost $r$-internal. Then there is a set $D \supset A$, some $b \models p$ and some realizations $c_1, \cdots , c_m \models r$ such that $b \forkindep_A D$ and $b \in \acl(c_1, \cdots , c_m , f(b), D)$. We pick $m$ minimal (among all such $m$ and $D$). Minimality of $m$ implies that $(c_1, \cdots , c_m) \models r^{(m)} \vert_D$, and with some forking calculus, that $c_1, \cdots , c_m \forkindep_D f(b)$.

    As a consequence, we see that $(c_1, \cdots , c_m) \models (r \vert_{f(b)D})^{(m)} = (r^{(m)})\vert_{f(b)D}$, and in particular $(c_1, \cdots , c_m , f(b)) \models ( r^{(m)} \otimes f(p))\vert_D$.

    We now show that the tuples $(c_1, \cdots, c_m, f(b))$ and $b$ are interalgebraic over $D$. We already know that $b \in \acl(c_1, \cdots , c_m, f(b),D)$. Suppose, for a contradiction, and without loss of generality, that $c_1 \not\in \acl(bD)$. This implies that $c_1 \forkindep_A bD$ by minimality. As $b \forkindep_A D$ we obtain $c_1 b \forkindep_A D$, and therefore $b \forkindep_A c_1 D$. This contradicts minimality of $m$, as we could then replace $D$ by $c_1 D$. So $c_i \in \acl(bD)$ for all $1 \leq i \leq m$.

    This interalgebraicity gives us:
    \begin{align*}
        n & = U(b/f(b)A) \\
        & = U(b/f(b)D) \\
        & = U(c_1, \cdots , c_m / f(b)D) \\
        & = m \text{ .}
    \end{align*}

    The implication $(2) \Rightarrow (3)$ is immediate. As for the implication $(3) \Rightarrow (1)$, a proof similar to that of Proposition \ref{pro: domeq for internal} gives that if $p \domequiv r^{(n)} \otimes f(p)$, then there are $b \models p$, a model $M \models T$ and some $c_1 , \cdots , c_n \models r$ such that $b$ and $c_1, \cdots ,c_n$ are interalgebraic over $f(b)M$ and $b \forkindep_A M$, which implies that $f$ is uniformly almost $r$-internal.

    Finally, the equivalence $(4) \Leftrightarrow (1)$ follows from \cite[Proposition 3.16]{jaoui2023relative} and is easy to prove.
\end{proof}

From Theorem \ref{theo: unif-int-dom-char}, it is easy to construct other examples where $m$ is $0$, but not because $p$ is orthogonal to the fibers. We first give an example, pointed out to us by the anonymous referee, which may be the simplest possible. It is \cite[Example 6.3.1, page 306]{Buechler_2017}:

\begin{example}
    Consider $G = (\mathbb{Z}/4\mathbb{Z})^{\omega}$, which is an $\aleph_1$-categorical group of Morley rank $2$. Let $p$ be its generic type, and consider the $\emptyset$-definable map $f : x \rightarrow 2x$, which gives us a map $f : p \rightarrow f(p)$. As proven in \cite[Example 6.3.1]{Buechler_2017}, the types $p$ and $f(p)$ are domination equivalent. Moreover, it is easy to see that the fibers of $f$ are internal to $f(p)$. So we have both $p \domequiv f(p)$ and $p \not\perp \fibre{p}{b}$ for any $f(b) \models f(p)$.
\end{example}

We also give a differential field example:

\begin{example}\label{ex: lod-der-of-constants}
    Work again in $\mathrm{DCF}_0$, and consider $p \in S(\mathbb{Q}^{\mathrm{alg}})$ the generic type of $\left( \frac{x'}{x}\right)' = 0$. If $r \in S(\mathbb{Q}^{\mathrm{alg}})$ is the generic type of the constant field $\calC$, then $p$ is the pullback of $r$ under the logarithmic derivative $\log_{\delta}$ in the sense of \cite{jin2020internality}, denoted $\log_{\delta}^{-1}(r)$. 

    The map $\log_{\delta} : p \rightarrow r$ is a fibration, and its fibers are minimal and $r$-internal. It is well-known (see, for example, \cite[Lemma 4.2]{chatzidakis2015differential}) that $p$ is not almost $r$-internal. By Theorem \ref{theo: unif-int-dom-char}, this implies that $p$ is not domination equivalent to $r \otimes \log_{\delta}(p)$, and therefore we must have $p \domequiv \log_{\delta}(p)$ (in other words $p \domequiv r$).
\end{example}

Example \ref{ex: lod-der-of-constants} can be used to construct fibrations $g : q \rightarrow g(q)$ with $q \domequiv r^{(m)} \otimes g(q)$ for some minimal $r$ and any $0 < m < n$, where $n$ is the $U$-rank of a fiber. 

Indeed, consider the type $p$ of Example \ref{ex: lod-der-of-constants} and $r$ the generic type of the constant field (all over $\mathbb{Q}^{\mathrm{alg}}$), and pick some integers $0 < m < n$. Consider the type $q = p^{(n-m)} \otimes r^{(m)}$ as well as the map $g$ given by:
\begin{align*}
    g: q & \rightarrow r^{(n-m)} \\
    (x_1, \cdots , x_{m-n}, y_1, \cdots , y_m) & \rightarrow (\log_{\delta}(x_1) , \cdots , \log_{\delta}(x_{n-m}))
\end{align*}
The type $q$ has $U$-rank $2n-m$, the map $g$ has $r$-internal fibers of $U$-rank $n$. Moreover, as the type $p$ is domination equivalent to $r$ by Example \ref{ex: lod-der-of-constants}, we have that $q \domequiv r^{(n-m)} \otimes r^{(m)} \domequiv r^{(m)} \otimes g(q)$ as $g(q)$ is the type of $n-m$ generic independent constants.

 \subsection{Domination decomposition and reduction}\label{sub: dom-and-red}

Recall that by Proposition \ref{pro: domeq for internal}, if $r \in S(A)$ is a minimal type and $p \in S(A)$, then it is almost $r$-internal if and only if $p \domequiv r^{\left( U(p) \right)}$. In this subsection, we prove a finer connection between internality and domination. 

First, using standard techniques (see for example \cite[Proposition 3.4.12]{wagner2000simple}), we obtain definable maps to internal types from non-domination:

\begin{lemma}\label{lem: map-from-no-dom-fix}
    Let $p\in S(A)$ be stationary and consider some $A$-definable map $f : p \rightarrow f(p)$. Suppose that there is some minimal type $r$ and some $n > 0$ such that $p \trianglerighteq r^{(n)}$ but $f(p) \not\trianglerighteq r^{(n)}$. Then there is an $A$-definable map $g : p \rightarrow g(p)$ such that:
    \begin{itemize}
        \item $g(p)$ is internal to the family $\calR$ of $A$-conjugates of $r$,
        \item for any $a \models p$, we have $g(a) \not\in \acl(f(a)A)$.
    \end{itemize}
\end{lemma}

\begin{proof}

    The assumption implies in particular that $f(p) \not\trianglerighteq p$, and that $r$ appears in the domination-decomposition of $p$ with a higher power than in the one of $f(p)$. By Lemma \ref{lem: min-dom} (and the remark following it), there is some $\aleph_1$-saturated $M \supset A$, some $a \models p\vert_M$ and some $e \models r \vert_M$ such that $e \forkindep_M f(a)$ and $e \not\forkindep_M a$. 

    Consider $\Cb(\stp(eM/aA))$, and let $d$ be a finite subtuple of maximal $U$-rank over $A$ (which exists as $\tp(a/A)$ has finite $U$-rank). Then $d$ is $\mathcal{R}$-internal. Indeed, we know that $d \in \dcl(e_1M_1, \cdots , e_nM_n)$ for some Morley sequence $e_1M_1, \cdots , e_nM_n$ in $\stp(eM/aA)$. An easy induction and forking computation shows that $a \forkindep_A M_1 \cdots M_n$, which, as $d \in \acl(aA)$, gives us $d \forkindep_A M_1 \cdots M_n$. This implies that $\tp(d/A)$ is $\calR$-internal, as the $e_i$ are realizations of $\calR$.

    Moreover, we know that $d \in \acl(aA)$ by properties of canonical bases. Pick $d'$ to be the canonical parameter of the finite set of realizations of $\tp(d/aA)$. Then $d' \in \dcl(aA)$ and is $\mathcal{R}$-internal. This gives a definable function $g : p \rightarrow g(p)$ to the internal type $g(p) = \tp(d'/A)$. To conclude the proof, we just have to show that $d' \not \in \acl(f(a)A)$. Note that $d \in \acl(d'A)$, so it is enough to show that $d \not \in \acl(f(a)A)$. Suppose on the contrary that $d \in \acl(f(a)A)$. As $\Cb(\stp(eM/aA)) \in \acl(dA)$, this implies that $eM \forkindep_{f(a)A} aA$, which gives us $e \forkindep_{f(a)M} a$, and as $e \forkindep_M f(a)$, that $e \forkindep_M a$, a contradiction.

\end{proof}

This allows us to connect domination to the \emph{reduction} of a type (see \cite{moosa2011model} and \cite{jin2018constructing}):

\begin{definition}
    Let $a$ be a tuple and $\calP$ be any $A$-invariant family of (potentially partial) types. A $\calP$-reduction of $a$ over $A$ is a (unique up to interalgebraicity) tuple $b \in \acl(aA)$ such that $\tp(b/A)$ is almost $\calP$-internal and for any other $b' \in \acl(aA)$, if $\tp(b'/A)$ is almost $\calP$-internal, then $b' \in \acl(bA)$.
\end{definition}

Remark that, up to interalgebraicity, we may as well assume that $\tp(b/A)$ is $\calP$-internal. As Jin notes in \cite{jin2018constructing}, if $\tp(a/A)$ has finite $U$-rank, then a $\calP$-reduction of $a$ over $A$ always exists: if $b \in \acl(aA)$ is a $\calP$-internal tuple of maximal $U$-rank, then $b$ is a $\calP$-reduction of $a$ over $A$. 

Moreover, we can find a reduction in $\dcl(aA)$ instead: if $b \in \acl(aA)$ is a $\calP$-reduction, let $b_1, \cdots , b_n$ be its orbit under $aA$-automorphisms, and let $b'$ be the canonical parameter of $\{ b_1, \cdots , b_n \}$. Then $b' \in \dcl(aA)$ and $\tp(b'/A)$ is $\calP$-internal. Since $b \in \acl(b')$, we also have $U(b'/A) \geq U(b/A)$, therefore $b'$ is a $\calP$-reduction of $a$ over $A$.

We prove:

\begin{theorem}\label{theo: dom-as-red}
   Let $\tp(a/A) = p$ be a non-algebraic stationary type, and assume that we have $p \domequiv r_1^{(n_1)} \otimes \cdots \otimes r_l^{(n_l)}$ where the $r_i$ are pairwise orthogonal minimal types (and all $n_i > 0$). For each $i$, let $\calR_i$ be the family of $A$-conjugates of $r_i$. Then $n_i$ is the $U$-rank of the $\calR_i$-reduction of $a$ over $A$. 
\end{theorem}

\begin{proof}
     Fix some $r_i$, and let $m_i$ be the $U$-rank of the $\calR_i$-reduction of $a$ over $A$. This is given by some $A$-definable map $\pi: p \rightarrow \pi(p)$ such that $\pi(p)$ is $\calR_i$-internal of maximal $U$-rank $m_i$. By a well-known fact (see \cite[Lemma 5.4]{freitag2025finite} for example), this implies that $\pi(p)$ is almost $r_i$-internal, and thus that $n_i \geq m_i$. We will now show the other inequality. 
     
     By \cite[Corollary 7.4.6]{pillay1996geometric}, there is an $A$-definable map $f : p \rightarrow f(p)$ such that $f(p)$ is $\calR_i$-internal. If $U(f(p)) \geq n_i$, we are done. Otherwise, by Lemma \ref{lem: map-from-no-dom-fix}, there is an $A$-definable map $g : p \rightarrow g(p)$ such that $g(a) \not\in \acl(f(a)A)$ and $g(p)$ is $\calR_i$-internal. Consider the map $f \times g$ sending any $b \models p$ to $(f(b),g(b))$. Its image is $\calR_i$-internal, and of $U$-rank strictly greater than $f(p)$. Repeating this process, we must eventually find some $A$-definable function $h : p \rightarrow h(p)$ such that $U(h(p)) \geq n_i$ and $h(p)$ is $\calR_i$-internal, proving $m_i \geq n_i$.
\end{proof}

While writing this article, we realized that this theorem is in fact a consequence of Buechler's work on levels. More precisely, Buechler defines the first level of some $b$ over $A$ as the set:
\[l_1(b/A) = \{ c \in \acl(bA) : \tp(c/A) \text{ is semiminimal} \}\]
\noindent and shows in \cite[Proposition 3.1]{buechler2008vaught} that $l_1(a/A)$ dominates $a$. This yields Theorem \ref{theo: dom-as-red} directly. We chose to leave this theorem and its applications because our methods are slightly different, and because we obtain a type-by-type version of his result, instead of considering the entire first level at once. Also note that Buechler's theory of levels has been generalized to simple theories by Palac\'{i}n and Wagner in \cite{palacin2013ample}, see their Definition 3.1, and their Theorem 3.6 for a generalization of Buechler's \cite[Proposition 3.1]{buechler2008vaught} in that context.

\begin{example}
    (\cite{jaoui2023relative}, Example 5.6). In DCF$_0$, let $p$ be the generic type over $\mathbb{C}$ of the system:
    \[\begin{cases}
        x'=x^2(x-1)\\
        y'=xy
    \end{cases}\] and let $f$ be the fibration given by projection on the $x$ coordinate. It is shown in \cite[Corollary 5.5]{jaoui2023relative} that $f:p\to f(p)$ is uniformly almost internal to $r$ the generic type type of the constants. By Theorem \ref{theo: unif-int-dom-char}, $p\domequiv r\otimes f(p)$. Note that since $U(p)=2$, $U(f(p))=1$ and hence is minimal. Also note that $f(p)\perp r$ by the results of \cite{rosenlicht1974nonminimality} or \cite[Example 2.20]{hrushovski2003model}. Hence for some (any) $a\models p,$ the $\calC$-reduction of $a$ over $\mathbb{C}$ and the $f(p)$-reduction of $a$ over $\mathbb{C}$ are both $U$-rank one. 
\end{example}

Theorem \ref{theo: dom-as-red} shows that finding the domination-decomposition should be a difficult task. However, we now point out how recent results of the authors \cite{eagles2024internality} can be used to find part of the domination decomposition for the generic types of some autonomous differential equations in $\mathrm{DCF}_0$. 

Fix some $\calU \models \mathrm{DCF}_0$, let $\calC$ be its field of constants, and $F$ be an algebraically closed subfield of $\calC$. Consider the following system of equations:

\[\begin{cases}
    y_1' = f_1(y_1, \cdots , y_k) \\
    \vdots \\
    y_k' = f_k(y_1, \cdots , y_k) 
\end{cases} \tag{$\dagger$} \label{sys of eq}\]
\noindent where $f_i \in F(x_1, \cdots , x_n)$ for all $i$, and let $p$ be its generic type. 

The data of such a system of equations is equivalent to that of a rational vector field $X$, i.e. a rational section of the tangent bundle of the affine space $\mathbb{A}^n$. Its \emph{Lie derivative} associates to any $g \in g(x_1, \cdots  x_n)$ the quantity:
\[\mathcal{L}_X(g) = \sum\limits_{i=1}^{k} \partials{g_j}{x_i}f_i\]
The results and methods of \cite{eagles2024internality} give:

\begin{theorem}\label{theo: red-for-autonomous}
    The Lascar rank of a $\calC$-reduction of $p$ is the maximal $k$ such that there exist rational functions $g_1 , \cdots , g_k \in F(x_1, \cdots , x_k) \in F(x_1, \cdots , x_k)$, algebraically independent over $F$, such that for all $1 \leq i \leq k$, either:
    \begin{enumerate}[A.]
        \item $\calL_X(g_i) = \lambda_i g_i$ for some $\lambda_i \in F$,
        \item $\calL_X(g_i) = 1$.
    \end{enumerate}
\end{theorem}

We only give a sketch of the proof, and direct the reader to \cite{eagles2024internality} for more details.

\begin{proof}[Proof sketch]
    The $g_i$ of the theorem give, by applying the chain rule, a map to the generic type of the system:

    \[\begin{cases}
        z_1 ' = \lambda_1 z_1 \text{ or } 1 \\
        \vdots \\
        z_{k}' = \lambda_{k} w_{k} \text{ or } 1 
    \end{cases}\]
\noindent where each line depends on whether we are in Case A. or B. The generic type of such a system is always $\calC$-internal of $U$-rank $k$. By similar methods as in the proof of \cite[Theorem 3.10]{eagles2024internality} (see also Lemma 3.1 in the same article), we see that if $f : p \rightarrow f(p)$ is an $F$-definable map to some $\calC$-internal type $f(p)$, then $f(p)$ has a definable map to the generic type of a system of that form. This yields the result.
\end{proof}

This theorem gives us the Lascar rank of the $\calC$-reduction of $p$. As a corollary of Theorem \ref{theo: dom-as-red}, we obtain that the generic type of the constants appears in the domination decomposition of $p$ with a Morley power of $k$.

\begin{example}
    Consider the classic Lotka-Volterra system:
    \[
    \begin{cases}
        x' = ax + b x y \\
        y' = ay + d x y
    \end{cases}
    \]
    for some $a,b,d \in \calC$ (so the coefficients of the linear terms are equal). It is shown in \cite[Theorem 4.7]{eagles2024internality} that the generic type $p$ of this system is 2-analyzable in the constants, but not almost internal to the constants. The $\calC$-reduction is given by the map $(x,y) \rightarrow \frac{x}{b} - \frac{y}{d}$, which goes to solutions of $z' = az$. Therefore $p \domequiv r$, where $r$ is the generic type of the constants.
\end{example}

\subsection{Dichotomy when the fibers have no proper fibrations}\label{sub: no-prop-fib}

Recall the following definition from \cite{moosa2014some}:

\begin{definition}
    A type $p = \tp(a/A)$ has no proper fibration if for any $b$, if $b \in \dcl(aA) \setminus \acl(A)$, then $a \in \acl(bA)$.
\end{definition}

It is proven in that article that any type without proper fibration is semiminimal. Moreover, it is not difficult to see that any finite rank type has a semiminimal analysis $f_1, \cdots , f_n$ such that the type of any fiber of an $f_i$ has no proper fibration. Therefore, it is natural to ask what our methods give when assuming that $f : p \rightarrow f(p)$ is a fibration such that the fibers $\fibre{p}{a}$, for $a \models p$, have no proper fibration. We obtain the following straightforward and expected consequence of Lemma \ref{lem: map-from-no-dom-fix}:

\begin{corollary}[Under assumption \ref{assump: non-orth}]
    Let $p \in S(A)$ be stationary and $f : p \rightarrow f(p)$ be a fibration such that its fibers have no proper fibration. Then either $p \domequiv f(p)$ or there is a minimal $r \in S(A)$ such that for any $f(a)$, the fiber $\fibre{p}{a}$ is almost $r$-internal and $p \domequiv r^{(n)} \otimes f(p)$, where $n = U(\fibre{p}{a})$.
\end{corollary}

\begin{proof}
    Fix some $a \models p$. The existence of some minimal $r_{f(a)}$ such that $\fibre{p}{a}$ is almost $r_{f(a)}$-internal is given by \cite[Proposition 2.3]{moosa2014some}. Assume that $p$ and $f(p)$ are not domination equivalent, then by Proposition \ref{prop: ortho-to-fib} and Corollary \ref{cor: suff-for-fib-descent}, there is some $r \in S(A)$ such that any fiber of $f$ is almost $r$-internal.

    Moreover, by Lemma \ref{lem: dom for internal fiber}, we have that $p \domequiv f(p) \otimes r^{(m)}$, for some $0 \leq m \leq n$, and since we assume that $p$ and $f(p)$ are not domination equivalent, we have $m > 0$. By Lemma \ref{lem: map-from-no-dom-fix}, there is an $A$-definable map $g : p \rightarrow g(p)$ such that $g(p)$ is $r$-internal and $g(a) \not\in \acl(f(a)A)$. Since $\fibre{p}{a}$ has no proper fibration, this implies that $a \in \acl(g(a)f(a)A)$.

    Since $g(p)$ is $r$-internal, there are $c_1, \cdots , c_k \models r$ and $B \supset A$ such that $g(a) \in \acl(c_1,\cdots , c_k, B)$ and $g(a) \forkindep_A B$. Pick some $B'$ such that $B' \equiv_{g(a)A} B$ and $B' \forkindep_{g(a)A} a$. Then $g(a) \in \acl(c_1', \cdots , c_k' ,B')$ for some $c_1' , \cdots , c_k ' \models r$, and $a \forkindep_A B'$. As $a \in \acl(g(a),f(a),A)$, we obtain than $a \in \acl(c_1' , \cdots , c_k' , f(a), B')$, and by the previous independence, the fibration $f : p \rightarrow f(p)$ is uniformly almost $r$-internal. Therefore, by Theorem \ref{theo: unif-int-dom-char}, we have $p \domequiv r^{(n)} \otimes f(p)$.
\end{proof}

\section{Splitting dichotomy for disintegrated fibers}\label{sec: disintegrated}

In this short section, we examine the specific case of a fibration $f : p \rightarrow f(p)$ such that its fibers are internal to some minimal type with disintegrated geometry. 

We will show that the fibers appear in the domination-decomposition of $p$ if and only if the map $f$ \emph{almost splits}, a notion studied in \cite{jin2020internality}, \cite{jaoui2023relative} and \cite{eagles2025splitting}:

\begin{definition}
    Let $p \in S(A)$ and $f : p \rightarrow f(p)$ be an $A$-definable fibration. It is split (resp. almost split) if there is a type $r \in S(A)$ such that for any $a \models p$, there is $b \models r$ such that $\dcl(aA) = \dcl(f(a)bA)$ (resp. $\acl$).
\end{definition}

Splitting can be seen to be equivalent to the existence of an $A$-definable map $p \rightarrow r \otimes f(p)$ such that the diagram
\begin{center}
\begin{tikzcd}
    p \arrow[rr] \arrow[rd, "f"] & &  r \otimes f(p) \arrow[dl] \\
    &  f(p) & 
\end{tikzcd}
\end{center}
\noindent commutes.

Our first lemma shows that if the minimal type is over the base parameters, the fibration almost splits (we thank the anonymous referee for suggesting this proof):

\begin{lemma}\label{lem: r-implies-product}
    Let $p \in S(A)$ and $f : p \rightarrow f(p)$ a fibration such that its fibers are of $U$-rank $n > 0$ and almost $r$-internal, where $r \in S(A)$ is minimal and disintegrated. Then for any $a \models p$, there are $c_1, \cdots , c_n \models r^{(n)}$ such that $a$ and $(f(a), c_1, \cdots ,c_n)$ are interalgebraic over $A$ and $f(a) \forkindep_A c_1, \cdots , c_n$
\end{lemma}

\begin{proof}\sloppy
    Fix some $a \models p$. Since $\fibre{p}{a}$ is almost $r$-internal and $r$ is disintegrated, it is in fact $r$-algebraic by Lemma \ref{lem: int-to-disin-implies-alg}. So there are $c_1, \cdots , c_m \models r$ such that $a \in \acl(f(a), c_1 , \cdots , c_m, A)$. We pick $m$ minimal, which implies that $f(a) \forkindep_A c_1, \cdots , c_m$. Since $n > 0$, we must also have $m > 0$.

    It is enough to show that $c_1 , \cdots , c_m \in \acl(a A)$. Indeed this implies that $a$ and $f(a), c_1, \cdots , c_m$ are interalgebraic over $f(a)A$, and therefore $n = m$ as $U(\fibre{p}{a}) = n$.
    
    Assume, for a contradiction, that some $c_i$ is not algebraic over $aA$. Without loss of generality, we may assume that each $c_1, \cdots , c_l$ is not algebraic over $aA$ and that each of $c_{l+1}, \cdots, c_m$ is, for some $1 \leq l \leq m$. Let $\overline{c} = (c_1, \cdots , c_l)$ and consider $d = \Cb(\stp(\overline{c}/aA))$. Since $\overline{c}$ is not algebraic over $aA$, we can pick an infinite Morley sequence $\overline{c} = \overline{c}_0, \overline{c}_1, \cdots$ in $\stp(\overline{c}/aA)$, and some $k$ such that $d \in \dcl(\overline{c}_0, \cdots , \overline{c}_k)$. This implies that $\overline{c}_{k+1} \forkindep_{\overline{c}_0 \cdots \overline{c}_k} aA$, as $d$ is also the canonical base of $\stp(\overline{c}_{k+1}/aA)$. 
    
    From this, we deduce $\overline{c}_{k+1} \not\forkindep_{A} \overline{c}_0 \cdots \overline{c}_k$. Indeed, otherwise we would obtain $a \forkindep_A \overline{c}_{k+1}$, and as $c_{l+1}, \cdots , c_m \in \acl(aA)$, that $a \forkindep_{Ac_{l+1} \cdots c_m} \overline{c}_{k+1}$. But $a \in \acl(f(a)\overline{c}_{k+1}c_{l+1} \cdots c_m A)$, so this implies $a \in \acl(f(a)c_{l+1}\cdots c_mA)$, contradicting minimality of $m$ (or that $a \not\in \acl(f(a)A)$ if $l = m$).

    This means that there must be some minimal $j \geq 0$ such that $\overline{c}_{j+1} \not\forkindep_{A} \overline{c}_0 \cdots \overline{c}_{j}$. By minimality of $j$, triviality of $r$ and indiscernibility of the sequence, there are $s,t$ such that $c_{j,s}$, the $s$-th coordinate of $\overline{c}_j$, and $c_{0,t}$, the $t$-th coordinate of $\overline{c}_{0}$, are interalgebraic over $A$. 

    Again by indiscernibility, this means that for all $i$, the $s$-th coordinates $c_{i,s}$ and $c_{0,s}$ are interalgebraic over $A$. Since $(\overline{c}_i)$ is a Morley sequence over $aA$, we obtain that $c_s = c_{0,s} \in \acl(aA)$, a contradiction.
    
\end{proof}

Almost splitting can be thought of as the strongest possible descent for the fibers. At the other end, the fibers could be pairwise orthogonal. We now show that, if the fibers are almost internal to disintegrated types, these are the only two options. We will use the following:

\begin{lemma}\label{lem: non-orth-fibers-implies-r}
    Let $p \in S(A)$ be stationary and $f : p \rightarrow f(p)$ be a fibration such that each fiber $\fibre{p}{a}$ is almost internal to a minimal disintegrated $r_{f(a)} \in S(f(a)A)$. Suppose that for any $(f(a),f(b)) \models f(p)^{(2)}$, we have $\fibre{p}{a} \not\perp \fibre{p}{b}$. Then there is a minimal, disintegrated $r \in S(A)$ such that $\fibre{p}{a}$ is almost $r$-internal (and thus $r$-algebraic) for all $f(a) \models f(p)$.
\end{lemma}

\begin{proof}
    Pick any $(a,b) \models p^{(2)}$. Since $\fibre{p}{a}$ and $\fibre{p}{b}$ are non-orthogonal and almost internal to $r_{f(a)}$ and $r_{f(b)}$, Corollary \ref{cor: min ortho int} implies that $r_{f(a)} \not\perp r_{f(b)}$, and as these are disintegrated, that $r_{f(a)}\vert_{f(b)} \not\perp^w r_{f(b)}\vert_{f(a)}$ by \cite[Lemma 16.2.11]{Baldwin_2017}. By Lemma \ref{lem: int-to-disin-implies-alg}, the types $\fibre{p}{a}\vert_{f(b)}$ and $\fibre{p}{b}\vert_{f(a)}$ are actually $r_{f(a)}$ and $r_{f(b)}$ algebraic, and therefore we obtain $\fibre{p}{a}\vert_{f(b)} \not\perp^w \fibre{p}{b}\vert_{f(a)}$.

    This implies that there is $b' \models \fibre{p}{b}\vert_{f(a)}$ such that $a \not\forkindep_{f(a)f(b)A} b'$ (note that $a \models \fibre{p}{a}\vert_{f(b)}$), which implies $a \not\forkindep_{f(a)A} b'$. Since $b' \forkindep_{f(b)A} f(a)$ and $f(a) \forkindep_A f(b)$, we get $b' \forkindep_A f(a)$. Therefore we have obtained $a \models \fibre{p}{a}$ and $b' \models p \vert_{f(a)}$ such that $a \not\forkindep_{f(a)A} b'$. This means that $p \not\perp \fibre{p}{a}$, and Lemma \ref{lem: min-descent-loc-mod} implies that there is a minimal $r \in S(A)$ such that $r_{f(a)} \domequiv r$. Since disintegration is preserved under domination-equivalence, see for example \cite[Chapter XVII, Theorem 2.6]{baldwin2017fundamentals}, the type $r$ is disintegrated.
    
   We may pick the same $r$ for all $f(a) \models f(p)$, and since $r_{f(a)}$ is almost $r$-internal and $\fibre{p}{a}$ is almost $r_{f(a)}$-internal, we see that $\fibre{p}{a}$ is almost $r$-internal, and thus $r$-algebraic by Lemma \ref{lem: int-to-disin-implies-alg}, for all $f(a) \models f(p)$.
\end{proof}

We finally obtain the main result of this subsection:

\begin{theorem}\label{theo: trivial-fibers-dicho}
    Let $p \in S(A)$ be stationary and $f : p \rightarrow f(p)$ be a fibration such that each fiber $\fibre{p}{a}$ is infinite and almost $r_{f(a)}$-internal, for some minimal disintegrated type $r_{f(a)} \in S(f(a)A)$. Then one of the following mutually exclusive possibilities must be true:

    \begin{enumerate}[(a)]
        \item for any $(f(a),f(b)) \models f(p)^{(2)}$, we have $\fibre{p}{a} \perp \fibre{p}{b}$,
        \item $f$ almost splits, and in particular there exists a minimal disintegrated type $r \in S(A)$ such that $p$ is interalgebraic with $f(p) \otimes r^{(n)}$, where $n$ is the $U$-rank of any fiber of $f$.
    \end{enumerate}
\end{theorem}

\begin{proof}
    Immediate consequence of Lemma \ref{lem: non-orth-fibers-implies-r} and Lemma \ref{lem: r-implies-product}.
\end{proof}

It is easy to find examples of (b) by simply taking Morley products. As for examples of (a), many can be found in the theory $\mathrm{DCF}_0$ of differentially closed fields of characteristic zero. For example, consider the second Painlev\'{e} equation:
\[P_{II}(\alpha) : y'' + ty + \alpha \]
\noindent for some $\alpha \in \calC$. If $\alpha$ is generic, then this isolates a type $p_{\alpha}$, which is known to be strongly minimal and disintegrated. By results of Freitag and Nagloo in \cite{freitag2025algebraic}, for any $\alpha, \beta$ generic independent constants, we have $p_{\alpha} \perp p_{\beta}$ (and they prove similar results for other Painlev\'{e} families).

\textbf{Acknowledgements.} The authors are thankful to Rahim Moosa for multiple conversations the subject of this article. They are also grateful to the anonymous referee for numerous improvements throughout the article, in particular for suggesting the current proof of Lemma \ref{lem: r-implies-product}.

\textbf{Funding.} The first author was supported by the Natural Sciences and Engineering Research Council of Canada (NSERC) [reference number CGS D - 588616 - 2024]
\bibliography{biblio}
\bibliographystyle{plain}

\end{document}